\documentclass{amsart}
\usepackage{amssymb,amsmath,amsfonts,epsfig,latexsym,texdraw}
\usepackage[all]{xy}

 \newtheorem{theorem}{Theorem}[section]
\newtheorem{definition}[theorem]{Definition}
\newtheorem{proposition}[theorem]{Proposition}
\newtheorem{lemma}[theorem]{Lemma}

\theoremstyle{definition}

\newtheorem{example}[theorem]{Example}

\def\Z{\ensuremath{\mathbb{Z}}}

\def\P{\ensuremath{\mathbb{P}}}

\def\C{\ensuremath{\mathbb{C}}}

\def\R{\ensuremath{\mathbb{R}}}

\def\cc{\ensuremath{\mathcal{C}}}
\def\cd{\ensuremath{\mathcal{D}}}
\def\cm{\ensuremath{\mathcal{M}}}
\def\cmbar{\ensuremath{\overline{\cm}}}
\def\ccbar{\ensuremath{\overline{\mathcal{C}}}}

\def\<{\ensuremath{\langle}}
\def\>{\ensuremath{\rangle}}

\def\PP{P\!P}

\DeclareMathOperator{\Pic}{Pic}

\DeclareMathOperator{\Span}{Span}

\DeclareMathOperator{\Supp}{Supp}
\DeclareMathOperator{\Star}{Star} \DeclareMathOperator{\Sym}{Sym}

\DeclareMathOperator{\Trop}{Trop}

\DeclareMathOperator{\MMW}{MMW}

\DeclareMathOperator{\dist}{dist}
\def\ol#1{{\overline{#1}}}

\def\cl#1{\operatorname{cl}(#1)}
\begin{document}

\title[Tropical Intersection Theory from Toric Varieties]{Tropical Intersection Theory from Toric Varieties}

\author[Katz]{Eric Katz}
\address{Department of Mathematics, University of Texas at Austin, Austin, TX 78712}
\email{eekatz@math.utexas.edu}

\begin{abstract}
We apply ideas from intersection theory on toric varieties to tropical intersection theory. 
We introduce mixed Minkowski weights on toric varieties which interpolate between equivariant and ordinary Chow cohomology classes on complete toric varieties.  These objects fit into the framework of tropical intersection theory developed by Allermann and Rau.   Standard facts about intersection theory on toric varieties are applied to show that the definitions of tropical intersection product on tropical cycles in $\R^n$ given by Allermann-Rau and Mikhalkin are equivalent.  We introduce an induced tropical intersection theory on subvarieties on a toric variety.  This  gives a conceptional proof that the intersection of tropical $\psi$-classes on $\cmbar_{0,n}$ used by Kerber and Markwig computes classical intersection numbers.
\end{abstract}

\maketitle

\section{Introduction}

Tropical Geometry is a rapidly developing field of mathematics.  In it, algebraic varieties are replaced by piecewise-linear objects called tropical varieties.  These tropical varieties capture a great deal of the nature of the original variety.  One of the original motivations for tropical geometry was that enumerative properties survive tropicalizations and, therefore, enumerative problems in algebraic geometry have piecewise-linear analogs that give the same answer.  To this end, it is desirable that tropical varieties have a good intersection theory.

There have been tropical intersection theories developed for tropical varieties in $\R^n$.  The definition of Mikhalkin in \cite{M}  uses a stable intersection in that two tropical cycles are perturbed so they intersect transversely.  The intersection number produced is independent of the perturbation.
Allermann and Rau developed a tropical intersection theory that is built out of the operation of intersecting tropical cycles with tropical Cartier divisors.  This is analogous to the construction of classical intersection theory in \cite{IT}.  

Tropical intersection theory has been applied in a number of contexts.  One such application is the intersection of $\psi$-classes on $\cmbar_{0,n}$ by Kerber and Markwig \cite{KM}.    They define tropical analogs of all the concepts involved and compute the tropical intersection numbers which are seen to be the same as the classical intersection numbers.

In this paper, we prove that the Allermann-Rau and the Mikhalkin tropical intersection theories are equivalent.   We also introduce an induced tropical intersection theory.  For $i:V\hookrightarrow X(\Delta)$, a subvariety of a toric variety satisfying certain properties, we show that the degree $\deg(i^*(d_1\cup\dots\cup d_k)\cap [V])$ for certain $d_i\in A^1(X)$ can be computed using tropical geometry. This induced tropical intersection theory is applied to show that the tropical $\psi$-class intersection numbers on $\cmbar_{0,n}$ agree with the classical ones.  We hope that this induced intersection theory sheds light on which other intersection problems can be computed tropically.

The method for the equivalence proof is to first prove that the intersection theories are equivalent for tropical cycles that are also fans.  In this case, the tropical intersection theory can be related to Chow cohomology on the toric variety described by the fan.  In the case that the toric variety is complete, the intersection theory is combinatorial.  By a result of Fulton and Sturmfels \cite{IToTV}, Chow cohomology classes are described by Minkowski weights which are exactly tropical cycles.  The cup product is described by the fan displacement rule given by Fulton and Sturmfels.    Equivariant Chow cohomology classes are described by piecewise polynomials by a result of Payne \cite{P}.  Tropical Cartier divisors which are piecewise-linear functions supported on tropical cycles can be interpreted as a 
hybrid object that we call {\em mixed Minkowski weights} which interpolate between equivariant and ordinary Chow cohomology.  We define a map $\kappa$ from mixed Minkowski weights to Chow cohomology.  It turns out that the associated Weil divisor of a tropical Cartier divisor $(c,f)$ is exactly the Chow cohomology class $\kappa(c,f)$.  The Mikhalkin definition of intersection product becomes the cup product of Chow classes by the fan displacement rule .  The equivalence of the two definitions of intersection theory follows from the identity   $(c_1\cup c_2)\cap [X(\Delta)]=c_1\cap (c_2\cap [X(\Delta)]).$
By a straightforward argument, the intersection theory on general tropical cycles in $\R^n$ can be reduced to the fan case.

For the result on tropical $\psi$-classes, we use an embedding of $\cmbar_{0,n}$ into a toric variety discovered by Kapranov \cite{K} and studied in detail by Gibney and Maclagan \cite{GM}.  The intersection problem on $\cmbar_{0,n}$ can be related to an intersection problem on the toric variety where it is then described by tropical intersection theory.

Our proof of the equivalence of tropical intersection theories which are both purely combinatorial makes use of non-combinatorial facts from the theory of toric varieties and intersection theory.  A combinatorial proof was given by Rau \cite{Rau}.

We would like to thank Sean Keel, Michael Kerber, Diane Maclagan, Hannah Markwig, Sam Payne, Johannes Rau, and David Speyer for valuable discussions.

\section{Toric Varieties and Tropical Geometry}

We review some notions from the theory of toric varieties.  A toric variety $X(\Delta)$ is determined by a rational fan $\Delta$ in a vector space $N_\R$.  The torus 
$T$ acting on $X$ has a one-parameter subgroup lattice canonically isomorphic to a lattice $N\subset N_\R$.  
The character lattice of $T$ is denoted by $M=N^\vee$.  The properties of the toric variety are reflected in those of the fan.  $X(\Delta)$ is complete if and only if $|\Delta|=N_\R$ while
$X(\Delta)$ is smooth if and only if $\Delta$ is unimodular, that is, if every cone in $\Delta$ is generated 
by a subset of an integral basis for $N$.   Any fan $\Delta$ can be refined to a unimodular fan 
$\widetilde{\Delta}$ which induces a resolution of singularities $X(\widetilde{\Delta})\rightarrow X(\Delta)$ \cite{F}.  
This resolution algorithm can be executed so that cones that are already unimodular do not need to be 
subdivided.  If $\tau$ is a cone in $\Delta$, the orbit closure $V(\tau)$ is a toric variety given by a fan 
$\Delta_\tau$ whose cones are $\overline{\sigma}=(\sigma+N_\tau)/N_\tau$ for $\sigma\supseteq\tau$
where $N_\tau$ is the linear span of $\tau$ in $N$.

Let $X=X(\Delta)$ be a complete $n$-dimensional toric variety with torus $T$.  Both the Chow cohomology, $A^*(X)$ and the equivariant Chow cohomology, $A_T^*(X)$ have combinatorial descriptions.  

The description of $A^*(X)$ is given by Fulton and Sturmfels \cite{IToTV}.   Let $\Delta^{(k)}$ denote the set of all cones in $\Delta$ of codimension $k$.  If $\tau\in\Delta^{(k+1)}$ is contained in a cone $\sigma\in\Delta^{(k)}$, let $v_{\sigma/\tau}\in N/N_\tau$ be the primitive generator of the ray $\overline{\sigma}$ in $\Delta_\tau$.

\begin{definition} A function $c:\Delta^{(k)}\rightarrow\Z$ is said to be a {\em Minkowski weight} if it satisfies the {\em balancing condition}, that is, for every $\tau\in\Delta^{(k+1)}$,
\[\sum_{\sigma\supset\tau} c(\sigma)v_{\sigma/\tau}=0\]
in $N/N_\tau$.
\end{definition}

The main result of \cite{IToTV} is that $A^k(X)$ is canonically isomorphic to the space of codimension $k$ Minkowski weights.  The usual operations on Chow cohomology can be described combinatorially. 
Let $c_1\in A^{k_1}(X),c_2\in A^{k_2}(X)$.  The Cartesian product $c_1\times c_2\in A^{k_1+k_2}(X(\Delta\times\Delta))$ is given by $(c_1\times c_2)(\sigma_1\times\sigma_2)=c_1(\sigma_1)c_2(\sigma_2)$.
The cup product is described by the fan displacement rule.  Let $v$ be a generic (as described in \cite{IToTV}) element of $N$.  
\[(c_1\cup c_2)(\gamma)=
\sum_{(\sigma_1,\sigma_2)\in\Delta^{(k_1)}\times\Delta^{(k_2)}}m_{\sigma_1,\sigma_2}^\gamma c_1(\sigma_1)c_2(\sigma_2)\]
where $m_{\sigma_1,\sigma_2}^\gamma$ are defined by
\[m_{\sigma_1,\sigma_2}^\gamma=
\begin{cases}
[N:N_{\sigma_1}+N_{\sigma_2}]&\text{if}\ \gamma\subset\sigma_1,\gamma\subset\sigma_2,\sigma_1\cap(\sigma_2+v)\neq \emptyset\\
0 &\text{otherwise}
\end{cases}
\]
The degree $\deg(c)$ of $c\in A^n(X(\Delta))$ is defined to be $c(0)$, the value of $c$ on the unique $0$-dimensional cone $0$.  

\begin{example}\label{fandisplacement}  Let $c_1,c_2$ be the codimension $1$ Minkowski weights below where the weight on each ray is $1$ except where indicated:
      
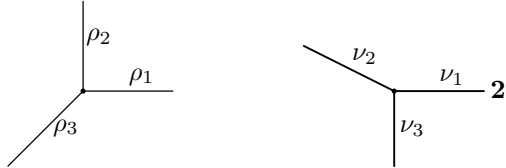
\begin{center}    
    \begin{texdraw}
       \drawdim cm  \relunitscale 0.4
       \linewd 0.05
        \move (3 3)        
        \fcir f:0 r:0.08
        \lvec (3 6)
        \move(3 3)
        \lvec(6 3)
        \move(3 3)
        \lvec(.5 .5)
        \move(4.5 3.2)\htext{$\rho_1$}
        \move(3.1 4.5)\htext{$\rho_2$}
        \move(2 1.5 )\htext{$\rho_3$} 
      \end{texdraw}     
    \hspace{1.5cm}
   \begin{texdraw}
       \drawdim cm  \relunitscale 0.4
         \linewd 0.06
        \move (3 3)        
        \fcir f:0 r:0.08
        \lvec (3 .5)
        \move(3 3)
        \lvec(6 3)
        \move(3 3)
        \lvec(0 4.5)
        \move(4.5 3.2)\htext{$\nu_1$}
        \move(6.2 2.75)\htext{${\bf 2}$}
        \move(1.6 3.9)\htext{$\nu_2$}
        \move(3.2 1.5 )\htext{$\nu_3$}  
       \end{texdraw}     
\end{center}
They are weights on some complete fan $\Delta$ in $\R^2$ that refines both of the fans indicated.  $\nu_2$ is in the direction $\<-2,1\>$  The Minkowski weight $c_1\cup c_2$ is determined by $(c_1\cup c_2)(0)$.  Let $v=\<2,1\>$.  The non-zero constants $m^0_{\rho\nu}$ correspond to the intersection of the first fan with the $\<2,1\>$-translate of the second:
      
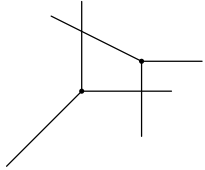
\begin{center}    
    \begin{texdraw}
       \drawdim cm  \relunitscale 0.4
       \linewd 0.05
        \move (3 3)        
        \fcir f:0 r:0.08
        \lvec (3 6)
        \move(3 3)
        \lvec(6 3)
        \move(3 3)
        \lvec(.5 .5)
        \move (5 4)        
        \fcir f:0 r:0.08
        \lvec (5 1.5)
        \move(5 4)
        \lvec(7 4)
        \move(5 4)
        \lvec(2 5.5)
       \end{texdraw}     
\end{center}
 Then $m^0_{\rho_1\nu_3}=[\Z^2:\<1,0\>+\<0,-1\>]=1$, $m^0_{\rho_2 \nu_2}=[\Z^2:\<-2,1\>+\<0,1\>]=2$.
 Therefore, $(c_1\cup c_2)(0)=m^0_{\rho_1\nu_3}c_1(\rho_1)c_2(\nu_3)+m^0_{\rho_2\nu_2}c_1(\rho_2)c_2(\nu_2)=1+2=3$.
\end{example}

Payne \cite{P} gave a combinatorial description of $A^*_T(X)$ even in the case that $X$ is not complete.  The ring of integral piecewise polynomial functions $\PP^*(\Delta)$ consists of continuous functions $f:|\Delta|\rightarrow\R$ such that the restriction $f_\sigma$ of $f$ to each cone $\sigma\in\Delta$ is a polynomial in $\Sym^*(M/(\sigma^\perp))$.
There is a canonical non-equivariant restriction map
\[\eta^*:A_T^*(X)\rightarrow A^*(X)\]
that is described combinatorially by Payne and the author \cite{KP}.  The kernel of that map contains $(M)$, the ideal generated by linear functions.  The case where $X$ is smooth was proven by Brion \cite{Brion}.
 Elements of $A^1_T(X)$ are given by $T$-Cartier divisors which are described by piecewise-linear functions on $\Delta$, that is
continuous functions $f:|\Delta|\rightarrow\R$ such that the restriction $f_\sigma$ of $f$ to each maximal cone $\sigma\in\Delta$ is an integer linear function in $M$. If $X$ is smooth, the canonical non-equivariant restriction map
\[\eta^*:A_T^1(X)\rightarrow A^1(X)\cong A_{n-1}(X)\]
has a much simpler combinatorial description as
\[\eta^*f=\sum_{\rho_i} -f(v_i)[V(\rho_i)]\]
where the sum is over rays $\rho_i$ in $\Delta$ and $v_i$ is the primitive integer vector along $\rho_i$ \cite{F}.

If $V(\tau)$ is an orbit closure in $X(\Delta)$ that is pointwise-fixed by the sub-torus $T(\tau)$, then there is a natural restriction map 
\[A_T^1(X)\rightarrow A_T^1(V(\tau))\]
A map $A^*_T(V(\tau))\rightarrow A^*_{T/T(\tau)}(V(\tau))$ is not canonically defined but requires a splitting of $T(\tau)\hookrightarrow T$.  Once a splitting is chosen, the image of a piecewise polynomial on $\Delta$ is a piecewise polynomial on $\Delta_\tau$, well-defined up to an element of $M(\tau)=\tau^\perp$.  The following diagram commutes
\[\xymatrix{
A^*_T(X)\ar[d]\ar[r]^{\eta^*}&A^*(X)\ar[d]\\
A^*_{T/T(\sigma)}(V(\tau))\ar[r]^>>>>>{\eta^*}&A^*(V(\tau)).}\]

Let $\Delta$ be a complete unimodular fan.  For $c\in A^k(X(\Delta))$, one can give a description of $\eta^*f\cup c\in A^{k+1}(X(\Delta))$.

\begin{lemma} \label{combformula} $\eta^*f\cup c\in A^{k+1}(X(\Delta))$ is the Minkowski weight given as follows, for $\tau\in\Delta^{(k+1)}$,
\[(\eta^*f\cup c)(\tau)=-\sum_{\sigma\in\Delta^{(k)}|\sigma\supset\tau} f(v_{\sigma/\tau})c(\sigma)\]
\end{lemma}

\begin{proof}
Let $i:V(\tau)\hookrightarrow X(\Delta)$ be the inclusion of the torus orbit.
After picking a splitting $T(\tau)\hookrightarrow T$, the lemma follows from the naturality of $\eta^*$ with respect to $i^*$.
\end{proof}

There is a connection between tropical varieties and Minkowski weights.  Let $V$ be a $l$-dimensional subvariety of a complete $X(\Delta)$ defined over $\C$.  $V$ is said to {\em intersect orbits properly} if and only if $\dim(V\cap V(\tau))=l+\dim(V(\tau))-n$ for all orbit closures $V(\tau)$.  The tropicalization of $V$, $\Trop(V)$ is an $l$-dimensional polyhedral fan in $N$ \cite{S}.  Moreover, the top dimensional cones $\sigma$ in $\Trop(V)$ are assigned a multiplicity $\omega(\sigma)$ that satisfies the balancing condition.
Proposition 2.2 of \cite{T} states that $V$ intersects orbits properly if and only if $\Trop(V)$ is a union of cones of $\Delta$.  Therefore, $\Trop(V)$ induces a Minkowski weight $c\in A^{n-l}(X)$ called the {\em associated cocycle}.  The class $c$ acts as an operational Poincar\'{e}-dual to $V$ in the following sense:
\begin{lemma}\cite[Lem 9.5]{Katz} Suppose $V$ intersects the torus orbits of $X$ properly. If $c$ is the associated cocycle of $V$, then
\[c\cap [X]=[V]\in A_l(X).\]
\end{lemma}

The following example will be useful in the sequel.
\begin{example} \label{diagonal} Let $\Gamma\subset(\C^*)^n\times (\C^*)^n$ be the diagonal.  Then a straightforward
 computation shows that $\Trop(\Gamma)$ is the diagonal in $\R^n \times \R^n$ with multiplicity $1$.  
 Therefore, if $\Delta$ is a complete rational fan in $\R^n\times\R^n$ so that the diagonal is a 
 union of cones, then the associated cocycle  of $\Gamma$ is a class $d\in A^n(X(\Delta))$ such that
\[d(\sigma)=\begin{cases}
1&\text{if $\sigma$ is contained in the diagonal}\\
0&\text{otherwise}
\end{cases}.\]
\end{example}

The diagonal is useful because if $\Delta'$ is a refinement of $\Delta\times\Delta$ for which $\Trop(\Gamma)$ is a union of cones and $\pi:X(\Delta')\rightarrow X(\Delta)$ is projection onto the first factor, then 
\[(c_1\cup c_2)\cap [X(\Delta)]=\pi_*((c_1\times c_2)\cap [\Gamma])=\pi_*((c_1\times c_2\cup d)\cap [X(\Delta')]).\]

A subvariety $V$ of $T$ can have additional properties that guarantee that some of the geometry of the closure of $V$ in $X(\Delta)$ can be seen from $\Trop(V)$. 

\begin{definition}\cite{T} The pair $(V,X(\Delta))$ is said to be {\em tropical} if the multiplication map $m:\overline{V}\times T\rightarrow X(\Delta)$ is faithfully flat and proper.  $(V,X(\Delta))$ is said to be {\em sch\"{o}n} if, in addition, $m$ is smooth.
\end{definition}

The fact that $(V,X(\Delta))$ is a tropical pair guarantees that $V$ intersects torus orbits properly.
It is proved in \cite{T} that if $V$ is sch\"{o}n for some tropical pair $(V,X(\Delta))$, it is sch\"{o}n for any tropical pair.    If in addition $V$ is sch\"{o}n and $X(\Delta)$ is smooth, then the stratification of $X(\Delta)$ as a toric variety pulls back to $V$ as a stratification with smooth strata.

\section{Mixed Minkowski Weights}

\begin{definition} Let $c \in A^k(X(\Delta))$ be a Minkowski weight of codimension $k$.  The {\em support} of $c$ is
\[\Supp(c)=\bigcup_{\ol{\sigma}|c(\sigma)\neq 0} \cl{\sigma},\]
the union of the closures of the cones on which $c$ is non-zero.
\end{definition}
Note that $\Supp(c)$ is a union of cones of codimension at least $k$.  The support of an associated cocycle of a subvariety $V$ of $X(\Delta)$ is $\Trop(V)$.  
 
\begin{definition} A mixed Minkowski weight $(c,f)$ of degree $k$ is a Minkowski weight, $c\in A^k(X(\Delta))$ together  with a continuous function $f:\Supp(c)\rightarrow\R$ that restricts to each cone as a linear function in $M/(\sigma^\perp)$.
\end{definition} 
Let $\MMW^k(\Delta)$ denote the set of degree $k$ mixed Minkowski weights.

\begin{lemma} If $\Delta$ is a unimodular complete fan, and $(c,f)$ is a mixed Minkowski weight then $f$ extends to a piecewise-linear function on $\Delta$.
\end{lemma}

\begin{proof}
Set $f$ to be zero on all rays not in $\Supp(c)$.  Extend $f$ linearly on top dimensional cones.  Such an extension agrees on overlaps.
\end{proof}

\begin{definition} For $\Delta$ unimodular, define the map
\[\kappa:\MMW^k(\Delta)\rightarrow A^{k+1}(X(\Delta))\]
as follows: for $(c,f)$, a mixed Minkowski weight, extend  $f$ to a piecewise-linear function on $\Delta$ and let $\kappa(c,f)=c\cup \eta^*f$.
\end{definition}

$\kappa(c,f)$ has a combinatorial description given by Lemma \ref{combformula}.  We can use this combinatorial formula to extend the definition to non-unimodular fans.  We will then show that the definition is independent of certain choices made.

\begin{definition} For complete $\Delta$, for $(c,f)\in\MMW^k(\Delta)$ and $\tau\in\Delta^{(k+1)}$, refine $\Delta_\tau$ to unimodular $\widetilde{\Delta}_\tau$, and let
\[\kappa(c,f)(\tau)=-\sum_{\rho\in\Delta_\tau^{(k)}} f(v_{\rho/0}) c(\Span^+(\tau,\rho))\]
\end{definition}

\begin{lemma} For $\Delta$ unimodular, $\kappa(c,f)$ is independent of the extension.
\end{lemma}

\begin{proof} Let $f_1,f_2\in\PP^1(\Delta)$ be extensions of $f$.
Then 
\[\kappa(c,f_1)(\tau)-\kappa(c,f_2)(\tau)=-\sum_{\sigma\in\Delta^{(k)}|\sigma\supset\tau}(f_1(v_{\sigma/\tau})-f_2(v_{\sigma/\tau})) c(\sigma).\]
For $\sigma$ with $c(\sigma)\neq 0$, we have $f_1|_\sigma=f_2|_\sigma$, so the sum is $0$.
\end{proof}

Now, we extend independence to the non-unimodular case.

\begin{lemma} $\kappa(c,f)(\tau)$ is independent of the choice of $\widetilde{\Delta}_\tau$.
\end{lemma}

\begin{proof}
Suppose we have two refinements $\widetilde{\Delta}_1,\widetilde{\Delta}_2$ of $\Delta_\tau$ with extensions $f_1,f_2$ of $\overline{f}$.   We may pick a common unimodular refinement $\widetilde{\Delta}$.  Then $f_1,f_2$ are both piecewise-linear on $\widetilde{\Delta}$ and the above lemma shows that $\kappa(c,f)(\tau)$ is independent of the choices involved.
\end{proof}

\section{Tropical Intersection Theory for Fans}

In this section, we prove that the Tropical Intersection Theory of \cite{AR} agrees with the fan displacement rule of \cite{IToTV}.  The main notions of the Tropical Intersection Theory of \cite{AR} have been introduced here in intersection theoretic guise.
A codimension $k$ tropical cycle supported on a fan $\Delta$ is just a codimension $k$ Minkowski weight on $\Delta$.  
A rational function on a tropical cycle $c$ is a piecewise-linear function $f$ on $\Supp(c)$ such that $(c,f)$ is a mixed Minkowski weight.  The associated Weil divisor of $f$ is manifestly the Minkowski weight $\kappa(c,f)$.

Given a codimension $k$ tropical cycle $c$ supported on $\Delta\subset N_\R$ and a morphism of fans $h:\Delta\rightarrow\Delta'$ induced from a homomorphism $h:N\rightarrow N'$ of lattices of dimension $n$ and $n'$, Gathmann, Kerber, and Markwig \cite{GKM} define a codimension $k+n'-n$ pushforward tropical cycle, $h_*(c)$.  If $\Delta$ and $\Delta'$ are unimodular, this pushforward cycle coincides with the image of $c$ under the pushforward in cohomology, $h_*:A^{k}(X(\Delta))\rightarrow A^{k+n'-n}(X(\Delta'))$ as explained combinatorially in \cite{IToTV}.

\begin{definition} Consider Minkowski weights $c_1\in A^{k_1}(X(\Delta))$, $c_2\in A_{k_2}(X(\Delta))$.   Let $x_i,y_i$ be coordinates on $N\times N$ corresponding to a fixed basis of $M$ applied to each factor.  Let $\chi_i=\min(0,y_i-x_i)$.  Let $\Delta'$ be a unimodular fan refining $\Delta\times\Delta$ such that each $\chi_i$ is piecewise-linear on $\Delta'$.   The {\em Allermann-Rau intersection product} is the Minkowski weight
\[c_1\bullet c_2=\pi_*(\kappa(\dots\kappa(\kappa((c_1\times c_2),\chi_1),\chi_2),\dots)\chi_n)))\]
where $\pi:\R^{2n}\rightarrow \R^n$ is projection onto the first factor.  
\end{definition}

By our description of $\kappa$, we can rewrite the above as
\[c_1\bullet c_2=\pi_*(\eta^*(\chi_1\cup\dots\cup\chi_n)\cup (c_1\times c_2)).\]

\begin{example} Let us revisit Example \ref{fandisplacement}.  Let $\Delta_1,\Delta_2$ be the complete fans whose rays are the support of $c_1,c_2$, respectively.  $c_1\times c_2$ is supported on the $9$ two-dimensional cones of $\Delta_1\times \Delta_2$ corresponding to the product of rays.  $c_1\bullet c_2=\pi_*(\kappa(\kappa(c_1\times c_2),\chi_1),\chi_2)$.  We refine $\Delta_1\times\Delta_2$ to $\widetilde{\Delta}$ to ensure that $\chi_1$ is piecewise-linear.   This involves subdividing $\rho_1\times\nu_1$ by adding the ray through $\omega_1=(\<1,0\>,\<1,0\>)$ and subdividing $\rho_3\times\nu_2$ by adding the ray through $\omega_2=(\<-2,-2\>,\<-2,1\>)$.  Let $d=\kappa(\chi_1,c_1\times c_2)$.  $d$ may be supported on $8$ rays, the $6$ rays of $\Delta_1\times\Delta_2$ and on $\omega_1,\omega_2$.  $(c_1\bullet c_2)(0)=\kappa(d,\chi_2)=\sum_\rho \chi_2(v_{\rho/0})d(\rho)$.  $\chi_2$ is zero on all of the rays except for 
\[\chi_2(\rho_2\times 0)=-1,\ \chi_2(0\times\nu_3)=-1\]
so it suffices to compute $d(\rho_2\times 0), d(0\times \nu_3)$.  For $d(\rho_2\times 0)$, we observe that $\rho_2$ is contained in $\rho_2\times\nu_1$, $\rho_2\times \nu_2$, and $\rho_2\times\nu_3$.  $\chi_1$ is zero on $\rho_2\times\nu_1,\rho_2\times\nu_3$ and
\[d(\rho_2\times 0)=-\chi_1(v_{\rho_2\times\nu_2/\rho_2\times 0})((c_1\times c_2)(\rho_2\times\nu_2))=-(-2)(1)=2.\]
Similarly,
\[d(0\times\nu_3)=-\chi_1(v_{\rho_1\times\nu_3/0\times\nu_3})((c_1\times c_2)(\rho_1\times\nu_3))=-(-1)(1)=1.\]
It follows that
\[\kappa(d,\chi_2)=-(\chi_2(\rho_2\times 0)d(\rho_2\times 0)+\chi_2(0\times\nu_3)d(0\times\nu_3))=(-1(2)+-1(1))=3.\]
This agrees with the computation using the fan displacement rule.  
\end{example}

\begin{lemma} Let $\Delta'$ be a refinement of $\Delta\times\Delta$ so that $\chi_1,\dots,\chi_n$ are piecewise-linear on $\Delta'$.  Then
\[\eta^*(\chi_1\cup\dots\cup\chi_n)=d\]
where $d$ is the associated cocycle of the diagonal from Example \ref{diagonal}.
\end{lemma}

\begin{proof}
We may suppose that $\Delta'$ is unimodular.
We first show that $\eta^*\chi_i$ is the associated cocycle of $V(w_i-z_i)$ where $w_k,z_k$ are coordinates on $(\C^*)^n\times (\C^*)^n$.  By the naturality of $\eta^*$ with respect to the toric projection $(\C^*)^n\times (\C^*)^n\rightarrow \C^*\times\C^*$ onto the $i$th factor of each $(\C^*)^n$,  it suffices to consider the case of $n=1$.  It is straightforward to verify that the weight of $\eta^*\chi_i$ on any ray contained in the diagonal is $1$.   Since the the varieties $V(w_i-z_i)$ intersect transversely, the class of $\eta^*(\chi_1\cup\dots\cup \chi_n)$ is exactly the Poincar\'{e}-dual of their intersection which is the diagonal.\end{proof}

\begin{theorem}  The Allermann-Rau intersection product is equal to the cup product on Minkowski weights: $c_1\bullet c_2=c_1\cup c_2$.
\end{theorem}

\begin{proof}
Let $[\Gamma]$ be the cycle-class corresponding to the diagonal.
On $X(\Delta')$, we have the following equation among cycle classes
\[(\eta^*(\chi_1\cup\dots\cup\chi_n)\cup (c_1\times c_2))\cap [X(\Delta')]=(c_1\times c_2)\cap [\Gamma].\]
This pushes forward by $\pi_*$ to the Poincar\'{e}-dual of $c_1\cup c_2$.
\end{proof}

\section{Tropical Intersection Theory}

In this section, we consider the more general tropical intersection theory of tropical polyhedral complexes in $\R^n$.  There are two definitions of tropical intersection product: one given by Mikhalkin in \cite{M} and one given by Allermann-Rau in \cite{AR}.  In this section, we prove that they are equal.  

A {\em polyhedral complex} $\cc$ in $N_\R$ is a finite collection of polyhedral cells in $N_\R$ that contains the faces of each of its polyhedra and such that the intersection of two polyhedra is a common face.  A polyhedral complex is said to be {\em rational} if all of its cells are rational with respect to the lattice $N$.  An integer affine linear function on $N$ is a function of the form $f(v)=\<v,m\>+a$ for $m\in M=N^\vee$ and $a\in\R$.

Let $\cc$ be an rational polyhedral complex in $\R^n$,  $\tau$, a cell in $\cc$, and $w$, a 
point in the relative interior of $\tau$. Define $\cd$ to be the set of all polyhedra in $\cc$ that contain $\tau$ as a 
face.  For $\sigma\in \cd$, let $C_\sigma$ be the cone
\[C_\sigma=\{v\in\R^n|w+\epsilon v\in \sigma \ \text{for some}\ \epsilon>0\}.\]
{\em The star of $\cc$ at $\tau$}, is the fan $\Star_\cc(\tau)$ made up of the union of the $C_\sigma$'s.  The minimal cone in $\Star_\cc(\tau)$ is $\overline{\tau}=\Span(\tau-w)$.  The star is independent of the choice of $w$.  

\begin{definition} A $k$-dimensional {\em tropical cycle} is a purely $k$-dimensional rational polyhedral complex together with a weight function on top-dimensional cells
\[\omega:\cc^{(k)}\rightarrow\Z\]
such that for any $(k-1)$-dimensional polyhedron $\tau$,
the weight induced by $\omega$ is  a tropical cycle on $\Star_\cc(\tau)$.
\end{definition}
The degree of a $0$-dimensional tropical cycle is the sum of the weights on its points.  

\begin{definition} Let $\cc_1$ and $\cc_2$ be tropical cycles of dimension $k$ and $l$ respectively.  Define the {\em Mikhalkin intersection product} $\cc_1\bullet\cc_2$ \cite{M} as follows: subdivide $\cc_1$ and $\cc_2$ such that their intersection $\cc_1\cap\cc_2$ is a subcomplex of each; for each $(k+l)-n$-dimensional cell $\tau$ of $\cc_1\cap\cc_2$; let $\Delta_1=\Star_{\cc_1}(\tau)$, $\Delta_2=\Star_{\cc_2}(\tau)$; pick a complete fan $\Delta$ that contains refinements of $\Delta_1$ and $\Delta_2$ as subcomplexes; $\Star_{\cc_1}(\tau)$, $\Star_{\cc_2}(\tau)$ induce Minkowski weights $c_1,c_2$ on $\Delta$; the weight on $\tau$ in $\cc_1\cdot\cc_2$ is $(c_1\cup c_2)(\overline\tau)$.  
\end{definition}

The Allermann-Rau intersection product \cite{AR} is  developed using tropical divisors. 

\begin{definition} A tropical Cartier divisor $\phi$ on $\cc$ is the data of $(U_i,\phi_i)$ where the $U_i$'s are open subcomplexes of $\cc$ such that $\cup U_i=\cc$ and $\phi_i$ is an integral piecewise-linear function on $U_i$ such that on each component of $U_i\cap U_j$, $\phi_i-\phi_j$ is an integral affine linear function.
\end{definition}
Note that if $\tau$ is any cone in $\cc$, a tropical Carter divisor $\phi$ induces a tropical Cartier divisor $\phi_{\Star(\tau)}$ on $\Star_\cc(\tau)$.  One takes some $U_i$ that contains $\tau$.  The restriction of $\phi_i$ to the cells containing $\tau$ induces a function $\phi_{\Star_\cc(\tau)}$ on $\Star_\cc(\tau)$.  One subtracts off a constant to ensure that the function is $0$ at the origin.  This function is well defined up to an global integer affine-linear function.

\begin{definition} The associated Weil divisor $\phi$ on $\cc$, denoted by $\phi\cdot \cc$ is a tropical cycle supported on the cones of positive codimension in $\cc$.  For $\tau\in\cc^{(k-1)}$, let $\omega_{\phi\cdot\cc}(\tau)=\kappa(\Star_\cc(\tau),\phi_{\Star_\cc(\tau)})(\overline{\tau})$,  the weight of $\overline{\tau}$ in the cycle given by taking the associated tropical Weil divisor of $\phi_{\Star_\cc(\tau)}$ on $\Star_\cc(\tau)$.
\end{definition}

The associated Weil divisor is easily seen to be a tropical cycle.

If $h:N\rightarrow N'$ is a map of lattices, there is a induced map of tropical cycles.  Let $\cc$ be a tropical cycle in  $N_\R$ and $\cc'$ is an integral polyhedral complex in $N'$ such that $h:\cc\rightarrow\cc'$ is a surjective map of polyhedral complexes of the same dimension.  Subdivide $\cc$ so that the image of each polyhedron in $\cc'$ is contained in a polyhedron.  If $\tau'$ is a top-dimensional polyhedron in $\cc'$, let $\tau_1,\dots,\tau_k$ be its top-dimensional pre-images in $\cc$.  The weight on $\tau'$ is 
\[\omega_{h_*(\cc)}(\tau')=\sum \omega({h_i}_*\tau_j)\]
where $\omega({h_i}_*\tau_j)$ is the weight on $\overline{\tau}_j$ as the image of the following map of tropical cycles on fans,
\[h_i:\Star_{\cc}(\tau_j)\rightarrow\Star_{\cc'}(\tau').\]

\begin{definition}  Let $\cc_1$, $\cc_2$ be tropical cycles in $\R^n$, the Allermann-Rau tropical intersection product \cite{AR} is defined to be
\[\cc_1\bullet \cc_2=\pi_*(\chi_1\cdot\ldots\cdot\chi_n \cdot (\cc_1\times \cc_2))\]
where $\pi:\R^n\times\R^n\rightarrow \R^n$ is projection onto the first factor.
\end{definition}

\begin{proposition} The Mikhalkin intersection product is equivalent to the Allermann-Rau intersection product.
\end{proposition}

\begin{proof}
Each intersection product is supported on cones in $\cc_1\cap \cc_2$.  For the Mikhalkin intersection product, this follows by definition.  For the Allermann-Rau intersection product, one notes that $\chi_1\cdot\ldots\cdot\chi_n\cdot(\cc_1\times\cc_2)$ is supported on $\Gamma\cap(\cc_1\times\cc_2)$ where $\Gamma$ is the diagonal in $\R^n\times\R^n$.  Such cells push forward to cells in $\cc_1\cap\cc_2$.

These definitions both reduce to computing intersection products on the stars of polyhedra in the intersection of two tropical cycles.  Therefore, the two definitions are equivalent.
\end{proof}

\section{Induced Intersection Theory}

Let $i:V\rightarrow X(\Delta)$ be a smooth $k$-dimensional subvariety of a complete smooth toric variety over $\C$.  Under certain conditions, intersection computations on $V$ can be performed using only the combinatorics of $\Trop(V)$.  We will apply this in the next section to the moduli space of rational curves.

Suppose $V$ intersects torus orbits properly.  Let $c\in A^k(X(\Delta))$ be the associated cocycle of $V$.  Then $\Supp(c)=\Trop(V)$. Since $c\cap [X(\Delta)]=[V]$, for $d\in A^{n-k}(X(\Delta))$, $i^*d\cap [V]=i^*d \cap (c\cap [X(\Delta)])=(i^*d\cup c)\cap [X(\Delta)]$.  Therefore, $\deg(i^*d\cap [V])$ depends only on the image of $d$ in $A^*(X(\Delta))/(\ker(\cup c))$.

\begin{lemma} \label{checkorbits} Let $d_1,d_2\in A^1(X(\Delta)$.  Suppose that for all $\tau\in\Delta^{(k+1)}$ with $\tau\subset\Supp(c)$ that $\deg((d_1\cup c)\cap [V(\tau)])=\deg((d_2\cup c)\cap [V(\tau)])$.  Then
$d_1=d_2$ in $A^*(X(\Delta))/(\ker(\cup c))$.
\end{lemma}

\begin{proof}  
It suffices to show $d_1\cup c=d_2\cup c$ as Minkowski weights.  Let $\tau\in\Delta^{(k+1)}$.  If $\tau\subset\Supp(c)$, then $(d_1\cup c)(\tau)=(d_2\cup c)(\tau)$.  Now suppose $\tau\not\subset \Supp(c)$.  Write $d_i\cap [V(\tau)]=\sum a_{i,\sigma}[V(\sigma)]$ where the sum is over $\sigma\in\Delta^{(k)}$ with $\sigma\subset\tau$.  Note that $c(\sigma)=0$ for all $\sigma\subset\tau$.  Then
\[(d_i\cup c)\cap [V(\tau)]=c\cap (d_i\cap [V(\tau)])=\sum a_{i,\sigma}c([V(\sigma)])=0.\]
\end{proof} 

\begin{definition} Let $d\in A^1(X(\Delta))$.  A piecewise-linear function $f$ on $\Trop(V)$ is said to {\em lift $d$ on $V$} if for all $\tau\in\Delta^{(k+1)}$, $\tau\subset\Trop(V)$,
\[\deg(i^*d\cap [i^{-1}V(\tau)])=\kappa(c,f)(\tau).\]
\end{definition}

In other words, if we extend $f$ to a piecewise-linear function $f$ on $\Delta$ then
$\deg(i^*d\cap [i^{-1}V(\tau)])=\deg((\eta^*f\cup c)\cap [V(\tau)]).$

\begin{theorem} \label{comparedegree} Let $V$ be a $k$-dimensional sch\"{o}n subvariety of smooth $X(\Delta)$ and suppose $(V,X(\Delta))$ is a tropical pair.   Let $d_1,\dots,d_k\in A^1(X(\Delta))$.  Let $f_1,\dots,f_k$ be piecewise-linear functions on $\Trop(V)$ lifting $d_1,\dots,d_k$.  Then
\[\deg(f_1\cdot\ldots\cdot f_k\cdot\Trop(V))=\deg((i^*d_1\cup\ldots i^*d_k)\cap [V]).\]
\end{theorem}

\begin{proof}
From $V$ sch\"{o}n, we know $V$ intersects the torus orbits of $X(\Delta)$ transversely and that one has the identity of cycle classes on $X(\Delta)$, $i_*[i^{-1}V(\tau)]=[V \cap V(\tau)]$.  Consequently if $f_i$ is extended as a piecewise-linear function on $\Delta$,
$\deg((\eta^*f_i\cup c)\cap [V(\tau)])=\deg((d_i\cup c)\cap [V(\tau)])$ for all $\tau\in\Delta^{(k+1)},\tau\subset\Supp(c)$.  Therefore, $\eta^*f_i=d_i$ in $A^*(X(\Delta))/(\ker(\cup c))$.
Therefore, 
\[\deg(f_1\cdot\ldots\cdot f_k\cdot\Trop(V))=\deg((\eta^*f_1\cup\ldots\cup \eta^*f_k\cup c)\cap [X(\Delta)])=\deg((d_1\cup\ldots d_k)\cap [V]).\]

\end{proof} 

\section{Intersection Theory on Moduli of Curves}

In \cite{KM}, Kerber and Markwig define tropical psi-classes as piecewise-linear functions on a model of tropical $\cmbar_{0,n}$.  They use the Allermann-Rau tropical intersection theory to compute the top-dimensional intersections of the $\psi$-classes.   The answer they obtain is equal to the classical answer.  We give a non-computational proof that these numbers are equal. 

$\cmbar_{0,n}$ is the moduli space of stable genus $0$ curves with $n$ marked points. See \cite{Manin} for an elaboration of the results cited below.
$\cmbar_{0,n}$ has a stratification by combinatorial types of dual graphs.  The big open set consists of smooth rational curves while boundary divisors consist of curves of arithmetic genus $0$ with two irreducible components.  The zero-dimensional strata correspond to curves with trivalent dual graphs while the one-dimensional strata correspond to dual graphs with one quadrivalent vertex and the remaining vertices trivalent.  For $k\in[n]$,  let $E_k$ be the set of one-dimensional strata for which the leaf-label $k$ is incident to the quadrivalent vertex.

$\cmbar_{0,n}$ carries natural cohomology classes $\psi_1,\dots,\psi_n\in A^1(\cmbar_{0,n})$.  Let $\ccbar_{0,n}$ be the universal curve over $\cmbar_{0,n}$.  Let $\omega$ be the relative dualizing sheaf of $\ccbar_{0,n}\rightarrow\cmbar_{0,n}$ which should be thought of as a relative cotangent bundle.  Let $\sigma_1,\dots,\sigma_n:\cmbar_{0,n}\rightarrow \cc_{0,n}$ be the sections of the universal curve given by the marked points.  Then $\psi_i=c_1(\sigma_i^*\omega)$.  It is straightforward to show that $\deg(\psi_k\cdot [\cmbar_{0,4}])=1$.  By the naturality of $\psi$-classes, one can show that if $S$ is the closure of a $1$-stratum of $\cmbar_{0,n}$ then
\[\psi_k\cap [S]=\left\{\begin{array}{lr}
1&S\in E_k\\
0&S\notin E_k.
\end{array}
\right.
\]

Kapranov proved that  $\cmbar_{0.n}$ is the Chow quotient of the Grassmannian, 
$G^0(2,n)//T$ where $G^0(2,n)$ is a particular dense open subset of the Grassmannian of $2$-planes in $\C^n$ and $T$ is the $n$-torus dilating the coordinates in $\C^n$ \cite{K}.  The Grassmannian, $G(2,n)$ is embedded by Pl\"{u}cker  coordinates into $\P^{\binom{n}{2}-1}$.  This embedding is equivariant with respect to the torus action.  $G^0(2,n)$ is the inverse image of the torus $(\C^*)^{\binom{n}{2}-1}$.  Therefore, the Pl\"{u}cker embedding induces an embedding into a Chow quotient, $i:\cmbar_{0,n}\hookrightarrow \P^{\binom{n}{2}-1}//T$.  Now, $\P^{{\binom{n}{2}-1}}//T$ is a $N=(\binom{n}{2}-n)$-dimensional toric variety which we will denote by $X(\Delta')$.  

The tropicalization of $\cm_{0,n}$ in $(\C^*)^N\subset X(\Delta')$ is well understood.  We summarize the description given by Kerber and Markwig in \cite{KM} which relies on the work of Speyer and Sturmfels \cite{SS}.  It it given in terms of combinatorial types of tropical rational curves.  An $n$-marked rational tropical curve is a metric tree with $n$ leaves labelled by integers in the set $[n]=\{1,\dots,n\}$.  All edges of the tree except the leaves are given a length.   Let the coordinates of $\R^{\binom{n}{2}}$ be indexed by two element subsets of $[n]$.  Denote the set of all such two element subsets by $I$. 
Define the following map
\[\begin{array}{rclcl}
\Theta&:&\R^n&\rightarrow& \R^{\binom{n}{2}}\\
&&a&\mapsto&(a_i+a_j)_{\{i,j\}\in I}
\end{array}
\]
The one-parameter subgroup lattice of $X(\Delta')$ canonically lies in the vector space $\R^{\binom{n}{2}}/\Theta(\R^n)$.
The image of a tropical curve $C$ in $\Trop(\cm_{0,n})\subset\R^{\binom{n}{2}}/\Theta(\R^n)$ is
the point whose coordinates are $\dist(\{i,j\})_{\{i,j\}\in I}$ where $\dist$ is the distance between the leaves in the metric graph $C$.  The image of all tropical curves is all of $\Trop(\cm_{0,n})$.   The image of tropical curves of a fixed combinatorial type gives a cone in the fan $\Trop(\cm_{0,n})$.  

We state some facts about the embedding that are proved using combinatorial techniques in a recent paper by Gibney and Maclagan \cite{GM}.  Their work draws on the study of the tropical Grassmannian by Speyer and Sturmfels \cite{SS} and the theory of tropical compactifications developed by Tevelev \cite{T}.   

\begin{proposition}
\ \\
\begin{enumerate}
\item  \cite[Prop 5.8]{GM} The tropicalization of the image of the moduli of curves, $\Trop(\cm_{0,n})$ is a union of cones in $\Delta'$.

\item \cite[Prop 5.4]{GM} If $\Delta$ is the union of cones in $\Delta'$ that support $\Trop(\cm_{0,n})$ then 
$X(\Delta)$ is a smooth toric variety.  

\item \cite[Lem 5.10]{GM} 
The embedding $i:\cmbar_{0,n}\rightarrow X(\Delta)$ induces an isomorphism 
\[i^*:\Pic(X(\Delta))\rightarrow \Pic(\cmbar_{0,n}).\]
\end{enumerate}

\end{proposition}

By the resolution of singularities algorithm for toric varieties as presented in \cite{F}, we may find a unimodular fan $\widetilde{\Delta}$ refining $\Delta'$ without subdividing any cones in $\Delta$.  Therefore, there is an open embedding $j:X(\Delta)\rightarrow X(\widetilde{\Delta})$.  Moreover, the composition $\cmbar_{0,n}\rightarrow X(\Delta)\rightarrow X(\widetilde{\Delta})$ is a regular embedding.  

The embedding $\cmbar_{0,n}\rightarrow X(\Delta)$ has a number of nice properties from the point of view of tropical compactifications.

\begin{proposition}\ \\
\begin{enumerate}
\item \cite{K} (see \cite[Thm 5.5]{T}) $(\cm_{0,n},X(\Delta))$ is a tropical pair.

\item \cite[Thm 1.11]{HKT} $\cm_{0,n}$ is sch\"{o}n.
\end{enumerate}
\end{proposition}

This implies by Proposition 2.3 of \cite{T} that $\cmbar_{0,n}$ intersects torus orbits in $X(\Delta')$ properly. Let $c$ be the associated cocycle of $\cmbar_{0,n}\subset X(\widetilde{\Delta})$ so that $c\cap [X(\widetilde{\Delta})]=[\cmbar_{0,n}]$.
The stratification of $X(\Delta)$ as a toric variety pulls back to the stratification of $\cmbar_{0,n}$.

We treat $\psi_i$ as elements of $A^1(X(\Delta))$.  By \cite[Section 2.3]{Brion} as explained in \cite{KP}, $\psi_i$ can be lifted to an element of $A_T^1(X(\Delta))$ which can be extended as a piecewise-linear function to an element of $A_T^1(X(\widetilde{\Delta}))$.  This gives a lift of $\psi_i$ to  $A^1(X(\widetilde{\Delta}))$.  

Mikhalkin introduced a tropical $\psi$-class in \cite{MRat}.  His definition corresponds to a Minkowski weight on the codimension $1$ cones in an embedding of $\Trop(\cm_{0,n})$ in a fan.  
Kerber and Markwig define a piecewise-linear function on $\Trop(\cm_{0,n})$ lifting $\psi_k$.  For  a subset $J\subset [n]$ of cardinality $1<|J|<n-1$, they define a vector $\nu_J\in\R^{\binom{n}{2}}$ by 
\[(\nu_J)_{\{i,j\}}=\begin{cases}
1,& \text{if}\ |\{i,j\}\cap J|=1\\
0,& \text{else}.
\end{cases}\]
$\nu_J$ is the ray in $\Trop(\cm_{0,n})$  corresponding to tropical curves with one finite internal edge with the marked points in $J$ attached one end-point and the marked points in $[n]\setminus J$ attached to the other.
For $k\in[n]$, the set $V_k$ is defined to be
\[V_k=\{\nu_J|k\notin J\ \text{and}\ |J|=2\}.\]
The piecewise-linear function $f_k$ is defined to be the linear extension on cones of $\Trop(\cm_{0,n})$  of 
\[\nu_J\mapsto 
\begin{cases} 1&\text{if}\  \nu_J\in V_k\\
0 &\text{else}.
\end{cases}\]

Since $f_k$ is a piecewise-linear function on $\Trop(\cm_{0,n})$, $(c,f_i)$ is a mixed Minkowski weight on $X(\Delta)$.  The associated Weil divisor of $f_k$ on $\Trop(\cm_{0,n})$ is determined by the following proposition.

\begin{proposition} \cite[Prop 3.5]{KM} If $\sigma$ is a codimension $1$ cone of $\cmbar_{0,n}$, then the weight $\omega_\sigma$ on $\sigma$ in the tropical Weil divisor of $f_k$ is given by 
\[\omega(\sigma)=\left\{\begin{array}{ll}
\binom{n-1}{2}&\text{$\sigma$ corresponds to a stratum in $E_k$}\\
0&\text{otherwise}.
\end{array}
\right.\]
\end{proposition}

This implies that $f_k$ lifts $\binom{n-1}{2}\psi_k$ on $\cm_{0,n}$.  Let $\Phi_k=\frac{1}{\binom{n}{2}}f_k$.  Theorem \ref{comparedegree} implies the following:

\begin{proposition}
\[\deg(\Phi_{k_1}\cdot\ldots\cdot\Phi_{k_l}\cdot\Trop(\cm_{0,n}))=\deg(\prod_l \psi_{k_l} \cap[\cmbar_{0,n}]).\]
\end{proposition}

This is the intersection of tropical $\psi$-classes computed in \cite{KM}.

\bibliographystyle{plain}

\end{document}